\input amstex
\input epsf
\documentstyle{amsppt}
\define\fact{!}
\redefine\cong{\equiv}
\define\per{\operatorname{per}}
\define\lamda{\lambda}
\define\perm{\operatorname{per}}
\define\suchthat{\mid}
\redefine\alpha{e}
\redefine\beta{l}
\refstyle{C}

\topmatter 

\title The number of terms in the permanent and the determinant of a generic
circulant matrix \endtitle
\rightheadtext{The determinant of a generic circulant matrix}
\author Hugh Thomas \endauthor
\date September 2002 \enddate

\address Fields Institute, 222 College Street, Toronto ON, M5T 3J1, Canada 
\endaddress

\email hthomas\@fields.utoronto.ca \endemail

\abstract
Let $A=(a_{ij})$ be the generic $n \times n$ circulant
matrix given by $a_{ij}=x_{i+j}$, with subscripts on $x$ interpreted 
mod $n$.  Define $d(n)$ (resp. $p(n)$) to be the number of terms in the
determinant (resp. permanent) of $A$.  
The function $p(n)$ 
is well-known and has several combinatorial
interpretations.  The function $d(n)$, on the other hand, has not been studied
previously.
We show that when $n$ is a prime power, $d(n)=p(n)$.  

\endabstract
\endtopmatter

\document

\head Introduction \endhead

A square matrix is said to be a circulant matrix if its rows are successive
cyclic permutations of the first row.  
Thus, the matrix $A=(a_{ij})$ with $a_{ij}=x_{i+j}$, subscripts on $x$
being interpreted mod $n$, is a generic circulant matrix.  

If we expand $\det(A)$, we obtain a polynomial in the $x_i$.  We define $d(n)$
to be the number of terms in this polynomial after like terms have been 
combined.  Similarly, we define $p(n)$ to be the number of terms in 
$\perm(A)$, the permanent of $A$.  

The function $p(n)$ was studied in Brualdi and Newman [1], where it was 
pointed out that the
main result of Hall 
[3] shows that $p(n)$ coincides with the number of solutions to
$$ \align y_1+2y_2+\dots+ny_{n} &\cong 0 \pmod n \\ 
y_1+\dots + y_n&=n \tag1\endalign$$
in non-negative integers.
Using this formulation, they showed by a generating function argument that
$$p(n)= \frac 1{n} \sum_{d|n} \phi(\frac nd)\binom{2d-1}{d}.\tag2$$

Setting $w_i = \sum_{j=i+1}^n y_j$, and rewriting (1)
in terms of the $w_i$, we see that $p(n)$ is the number of non-increasing 
$(n-1)$-tuples
$(w_1,\dots,w_{n-1})$ with $n \geq w_1 \geq \dots \geq w_{n-1} \geq 0$, such
that $\sum_i w_i \cong 0 \pmod n$.  If we let $L$ be the $(n-1)$-dimensional
lattice consisting of $(n-1)$-tuples of integers whose sum is divisible by $n$,
then this expression for $p(n)$ amounts to the number of points from $L$ lying
in the simplex with vertices  $(0,0,\dots,0)$, $(n,0,0,\dots,0)$,
$(n,n,0,\dots,0)$, $\dots$, $(n,n,\dots,n)$.  

There is another combinatorial interpretation for $p(n)$, as follows. 
Consider  all possible necklaces consisting of $n$ white beads
and $n$ black beads,
where two necklaces are considered equivalent if they differ by a cyclic 
permutation.  A straightforward P\'olya counting argument shows that 
the number of such necklaces is given by the right-hand side of (2), and thus
that the number of necklaces equals the number of terms in $\perm(A)$,
though no explicit bijection is known.

In contrast, little is known about $d(n)$, despite the fact that its definition
seems at least as natural.  It is clear that $d(n) \leq p(n)$ since 
every term which appears in the determinant also appears in the permanent. 
However, some terms from the permanent could be absent from the determinant 
due to cancellation. 
In this paper
we establish the following theorem:

\proclaim{Theorem} If $n$ is a prime power, $d(n)=p(n)$.  \endproclaim

The following table of values for $d(n)$ and $p(n)$ suggests
that the converse may also be true.  This is still open.

$$\matrix n & d(n) & p(n) & & & n & d(n) & p(n) \\
          1 & 1    & 1    & & & 7 & 246 & 246 \\
          2 & 2    & 2    & & & 8 & 810 & 810 \\
          3 & 4    & 4    & & & 9 & 2704 & 2704 \\
          4 & 10 & 10     & & & 10 & 7492 & 9252 \\
          5 & 26 & 26     & & & 11 & 32066 & 32066 \\
          6 & 68 & 80     & & & 12 & 86500 & 112720 \endmatrix $$

\head Background on Symmetric Functions \endhead

The proof of the theorem uses the theory of symmetric functions.  In this 
section, we review the necessary definitions and results.  For more 
detail on symmetric functions, see Stanley [4]. 

Symmetric functions are power series (in our case, over $\Bbb{Q}$) 
in an infinite number of variables
$z_1,z_2,\dots$, such that for any $(b_1,\dots,b_k) \in \Bbb{N}^k$ 
the coefficient of $z_{i_1}^{b_1}\dots z_{i_k}^{b_k}$ does not depend on the
choice of distinct natural numbers $i_1,\dots, i_k$. 

We write $\lamda \vdash q$ to signify that $\lamda$ is a partition of 
$q$.  
The symmetric functions of degree $q$ form a vector space whose dimension
is the number of partitions of $q$.  There are several standard bases for 
them.  
We shall
need two here: $\{m_\lamda\suchthat \lamda \vdash q\}$  
and $\{p_\lamda\suchthat \lamda \vdash q\}$.  For 
$\lamda=(\lamda_1,\dots,\lamda_k)$ with $\lamda_1\geq\lamda_2\geq \dots\geq
\lamda_k$, $m_\lamda$ is the power series in which
$z_{i_1}^{\lamda_1}\dots z_{i_k}^{\lamda_k}$ occurs with coefficient one
for every distinct sequence of natural numbers $i_1,\dots,i_k$ and no other
terms occur.  The symmetric function $p_i$ is defined to be 
$z_1^i+z_2^i+\dots$, and the symmetric function $p_\lamda$ is defined to be 
$ p_{\lamda_1}\dots
p_{\lamda_k}$.  

As we remarked,  
$\{p_\lamda\mid \lamda \vdash q\}$ and $\{m_\lamda\mid \lamda \vdash q\}$ form bases for the symmetric functions
of degree $q$.  We will need to convert from the $m_\lamda$ basis to the
$p_\lamda$ basis, which we shall do using a result of 
E\u gecio\u glu and Remmel  [2], 
for which we must introduce some notation.  

For $\lamda \vdash q$, let $k(\lamda)$ denote the number of parts 
of $\lamda$.  If $\lamda=\langle1^{\beta_1}2^{\beta_2}\dots q^{\beta_q}
\rangle$, let $\lamda!=1!^{\beta_1}2!^{\beta_2}\dots q!^{\beta_q}$, and
let 
$z_\lamda=\beta_1\fact\dots \beta_q\fact
1^{\beta_1}2^{\beta_2}\dots q^{\beta_q}$.

The Ferrers diagram of $\lamda$ consists of a 
left-justified column of rows of boxes
of lengths $\lamda_1, \lamda_2, \dots, \lamda_k$.  For $\mu$ another partition
of $q$,  
a filling of $\lamda$ by $\mu$ is a way to cover the Ferrers diagram of
$\lamda$ in a non-overlapping manner by ``bricks'' consisting of 
horizontal sequences of boxes of length $\mu_1, \mu_2,\dots$.
The weight of a filling is the product over all rows of the length of the 
final 
brick in each row.  In deciding whether two fillings are the same, 
bricks of the same size are considered to be indistinguishable.
Thus, the four fillings of $\lamda=(4,2)$ by $\mu=(2,2,1,1)$,
having weights respectively 4, 2, 2, and 2, are:

\smallskip

$$\epsfbox{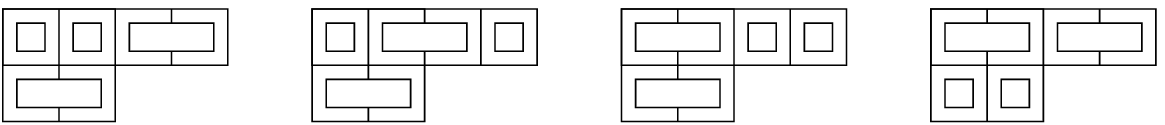}$$
 
\smallskip

Let $w(\lamda,\mu)$ be the sum over all distinct fillings
of $\lamda$ by $\mu$ of the weight of the filling. 
Then the result we shall need from [2] is that:

$$ m_\mu = \sum_{\lamda \vdash q} \frac{(-1)^{k(\mu)-k(\lamda)}w(\lamda,\mu)}
{z_\lamda}p_\lamda.$$

\head Proof of Theorem \endhead

First, let us consider what terms occur in $\per(A)$.  Let $b=(b_1,\dots,b_n)$
be an $n$-tuple of natural numbers summing to $n$.  Let $x^b$ denote 
$x_1^{b_1}x_2^{b_2}\dots x_n^{b_n}$.  We are interested in whether 
$x^b$ appears in $\per(A)$ with non-zero coefficient.  Suppose it does. 
Then there is some permutation 
$\sigma$ of $\{1,\dots,n\}$ such that $\prod_i x_{i-\sigma(i)}= x^b$.  
Since $\sum_i (i-\sigma(i))=0$, it follows that $\sum_i ib_i \cong 0
\pmod n$.  By the result of [3] mentioned above, this necessary condition
for $x^b$ to occur in $\per(A)$ is also sufficient.  

Now we proceed to consider $\det(A)$.  
Let $\xi$ be a primitive $n$-th root of unity.  
$A$ is diagonalizable, and its eigenvalues are $c_1,\dots,c_{n}$,
where $c_i=\sum_j \xi^{ij}x_j$.  Thus, $\det(A)=\prod c_i$.

Again, let $b=(b_1,\dots,b_n)$ be
an $n$-tuple of integers summing to $n$.
Let $q=b_1+2b_2+\dots+nb_n$.  
Let $\mu   =\langle 1^{b_1}2^{b_2}\dots n^{b_n}\rangle$, so 
$\mu \vdash q$.
Then $[x^b]\det(A)$, the coefficient
of  $x^b$ in $\det(A)$, is given by 
$$ [x^b]\det(A)=\sum_{f:\{1,\dots,n\}\rightarrow \{1,\dots,n\} \atop 
|f^{-1}(i)|=b_i} \xi^{\sum_{i=1}^n if(i)}.$$ 
We observe that this also equals $m_\mu(\xi^1,\xi^2,\dots,
\xi^{n},0,0,\dots)$.  (A symmetric function is a power series, so it is
not generally legitimate to substitute in values for the indeterminates,
but since $m_\mu$ is homogeneous and all but finitely many of the 
values being substituted are zero, it is allowed in this case.)  

Now we use the result of [2] mentioned above, which shows that:

$$[x^b]\det(A)=\sum_{\lamda \vdash q}\frac{(-1)^{k(\mu)-k(\lamda)}
w(\lamda,\mu)p_\lamda(\xi^1,\xi^2,\dots,\xi^{n},0,0,\dots)}{z_\lamda}.\tag3$$

What can we say about $p_\lamda(\xi^1,\xi^2,\dots,\xi^{n},0,0,\dots)$?
Firstly, $p_i(\xi^1,\xi^2,\dots,\xi^{n}$, $0,0,\dots) = n$ if $n|i$ and 
$0$ otherwise.  Thus, $p_\lamda( \xi^1,\xi^2,\dots,\xi^{n},0,0,\dots)=
n^{k(\lamda)}$ if all the parts of $\lamda$ are multiples of $n$, and $0$ 
otherwse.

(From this we could deduce  that if $ [x^b]\det(A)$
is non-zero, $q= b_1+2b_2+\dots+nb_n$ must be a multiple of $n$.  
Of course, we already know this, by the argument given when discussing
$\perm(A)$.)  

To establish our result, we must now show that if $q$ is a multiple of
$n$ and $n$ is a prime power, then the sum in (3) is non-zero. 
So assume that $n=p^r$, $p$ a prime.  
Let $v: \Bbb{Q}^\times \rightarrow \Bbb{Z}$ denote the usual valuation
with respect to $p$.

For $\lamda$ any partition of $q$, divide the fillings of $\lamda$ by $\mu$ 
into equivalence classes where
two fillings are equivalent if one can be obtained
from the other by rearranging the bricks within  each row, and by 
swapping the sets of bricks filling  pairs of rows of equal length.  

We wish to show that the contribution to (3) from the partition 
$\langle q \rangle$ (all of whose fillings form a single equivalence class)
has a smaller valuation than the sum of weights of the fillings in 
any equivalence class of fillings of any other partition.  Once this is
established, it follows that the sum (3) is non-zero.  

Fix $\lamda=(\lamda_1,\lamda_2,\dots)=\langle 1^{\beta_1}2^{\beta_2}\dots 
q^{\beta_q}\rangle $, with all the $\lamda_i$ divisible by 
$n$, and fix an equivalence class $\Cal{F}$ of
fillings of $\lamda$ by $\mu$.  Write $k$ for $k(\lamda)$.
Consider first a subclass of fillings 
$\Cal{G}$, those which can be obtained from some fixed $F \in \Cal{F}$ 
by rearranging the bricks in each row, but without interchanging rows.  
Let the $j$-th row of the filling $F$ of $\lamda$ contain $r_j$ bricks.  
Let $\alpha_{ij}$ be the number of these bricks having length $i$.   
The number of rearrangements of this row 
with ending in a brick of length $i$ is 
$\binom{r_j-1}{\alpha_{1j},\dots,\alpha_{ij}-1,\dots,\alpha_{nj}}$.

Thus, the total weight of all the rearrangements of this row is:

$$\align
\sum_{i=1}^n \binom{r_j-1}{\alpha_{1j},\dots,\alpha_{ij}-1,\dots,
\alpha_{nj}}i&=
\sum_{i=1}^n \frac1{r_j} \binom 
{r_j}{\alpha_{1j},\dots,\alpha_{ij},\dots,\alpha_{nj}}i\alpha_{ij}\\&=
\frac 1{r_j} \binom {r_j}{\alpha_{1j},\dots,\alpha_{nj}} 
\lamda_j.\endalign$$ 

It follows that the total weight of all the fillings in $\Cal{G}$ is:

$$
\prod_{j=1}^{k} \frac1{r_j} \binom {r_j}{\alpha_{1j},\dots,\alpha_{nj}} 
\lamda_j.$$

Consider the $\beta_i$ parts of $\lamda$ of length $i$.  The 
equivalence class $\Cal{F}$ determines a partition $\gamma(\Cal{F},i)$ of 
$\beta_i$, where the parts of $\gamma(\Cal{F},i)$ are the sizes of the
sets of rows of length $i$ which are filled with an indistinguishable set of
bricks.  
The sum of all the weights
over all the fillings in $\Cal{F}$ is:

$$\prod_{i=1}^q \frac{\beta_i!}{\gamma(F,i)!}
\prod_{j=1}^{k}\frac 1{r_j}\binom{r_j}{\alpha_{1j},\dots,\alpha_{nj}}
\lamda_j.$$

Writing $\delta(\Cal{F})$ for the partition of $k$ which is the 
sum of the $\gamma(\Cal{F},i)$ for all $i$,
the contribution to (3) from all these fillings is:

$$ \frac{(-1)^{k(\mu)-k}n^{k}}
{\delta(\Cal{F})!}
\prod_{j=1}^{k} \frac {(r_j-1)\fact}{\alpha_{1j}\fact
\dots\alpha_{nj}\fact} \tag4$$

Now let us consider (4) in the case where $\lamda=\langle q \rangle$.  
Here, we obtain

$$\frac{(-1)^{k(\mu)-1}n!}{b_1!\dots b_n!}.\tag5$$

We wish to show that (4) evaluated for any other equivalence class 
has a greater valuation with respect to $p$ than does (5).  
This is equivalent to showing that, for any $\lamda \ne \langle q \rangle$,
and any equivalence class of fillings $\Cal{F}$, that the following
expression has a positive valuation:

$$\align &
\frac{n^{k}}{\delta(\Cal{F})!}
\left(\prod_{j=1}^{k} \frac {(r_j-1)\fact}{\alpha_{1j}\fact
\dots\alpha_{nj}\fact}\right) \frac {b_1!\dots b_n!}{n!}\\
&\qquad = 
\left(\frac 1 {\delta(\Cal{F})!
} \prod_{i=1}^n \binom{b_i}{\alpha_{i1},\dots,
\alpha_{ik}}\right)\left(\frac{n^{k-1}}{(n-1)\dots(n-k+1)
\binom{n-k}{r_1-1,\dots,r_k-1}}\right).\tag6\endalign$$

We have written (6) as a product of two terms.  We will show that the
first term is an integer, and therefore has non-negative valuation, and
that the second term has positive valuation, which will complete the 
proof.  

For the first term, observe that if we rewrite the partition 
$b_i=\alpha_{i1}+\dots+\alpha_{ik}$ as $\langle 1^{c_{i1}}\dots n^{c_{in}}
\rangle$, 
then $c_{i1}!\dots c_{in}!$ divides $\binom{b_i}{\alpha_{i1},
\dots,\alpha_{ik}}$ because 

$$\frac {1}{c_{i1}!\dots c_{in}!}\binom{b_i}{\alpha_{i1},
\dots,\alpha_{ik}}$$
counts the number of ways of dividing $b_i$ objects into 
subsets of 
certain sizes, where we don't distinguish the different
subsets of the same size.  The first term of (6) is now disposed of by 
remarking that
$$\delta(\Cal{F})! | \prod_{i=1}^n c_{i1}!\dots c_{in}!$$
since each term in $\delta(\Cal{F})!$ implies the existence of at least one 
equal term among the product of factorials on the right hand side.  

For the second term, we need the following simple lemma:

\proclaim{Lemma} If $m<p^s$, 

(i) $v(\binom m d) < s$,

(ii) $v(\binom m {d_1,\dots,d_k})< (k-1)s$. \endproclaim

\demo{Proof} $$ \align
v(\binom m d)&= \left(\lfloor \frac m p \rfloor - \lfloor \frac
d p \rfloor - \lfloor \frac {(m-d)} p\rfloor \right) + \dots + 
\left(\lfloor \frac m {p^{s-1}} \rfloor - \lfloor \frac
d {p^{s-1}} \rfloor - \lfloor \frac {(m-d)} {p^{s-1}}\rfloor \right)\\
&\leq 1 + \dots + 1 = s-1. \endalign$$
The second part follows immediately from repeated application of the first
part, which completes the proof of the lemma.  \enddemo

Now, we know that $$\align v(n^{k-1})&=(k-1)r,\\ 
v(\binom{n-k}{r_1-1,\dots,r_k-1}) &\leq (k-1)(r-1),\\  
v((n-1)\dots(n-k+1))&=v((k-1)!)=\lfloor k/p\rfloor +\dots 
<\frac{k-1}{p-1} \leq k-1,\endalign$$
and the desired result follows.   

\head Acknowledgements \endhead

I would like to thank Graham Denham for a useful suggestions and 
Jim Propp for asking the question and for 
helpful comments on a previous draft.


\Refs

\widestnumber\key{4}

\ref
\key 1
\by R. Brualdi and M. Newman
\paper An enumeration problem for a congruence equation
\jour J. Res. Nat. Bur. Standards Sect. B 
\vol 74B \yr 1970 \pages 37--40
\endref

\ref
\key 2
\by \"O. E\u gecio\u glu and J. Remmel 
\paper Brick tabloids and the connection matrices between bases of
symmetric functions 
\jour Discrete Appl. Math.
\vol 34 \yr 1991 \pages 107--120
\endref

\ref
\key 3
\by M. Hall, Jr.
\paper A combinatorial problem on abelian groups
\jour Proc. Amer. Math. Soc. \vol 3 \yr 1952 \pages 584--587 
\endref

\ref
\key 4
\by R. Stanley
\book Enumerative Combinatorics
\bookinfo Volume 2
\publ Cambridge University Press
\publaddr Cambridge
\yr 1999
\endref

\endRefs
\enddocument